\documentclass[12pt,a4paper]{amsart} 

\usepackage{amssymb} 
\usepackage{amsmath} 
\usepackage{amsfonts}

\usepackage{epsfig}

\theoremstyle{plain} 
\newtheorem{theo} {Theorem} [section]

\newtheorem{lemm} [theo]{Lemma} 
 
\newtheorem{ques} [theo]{Question} 

\def\Int{\mathop{\rm Int} \nolimits} 
 
\def\RR{\mathbb{R}} 
\def\NN{\mathbb{N}} 
\def\ZZ{\mathbb{Z}}

\begin{document} 

\bigskip

\title{Embeddability of multiple cones} 
\author{D.~Repov\v s} 
\address{Institute of Mathematics, Physics and Mechanics, University of Ljubljana,  P.O.Box 2964, Ljubljana 1001, Slovenia} 
\email{dusan.repovs@guest.arnes.si} 
\author{W.~Rosicki} 
\address{Institute of Mathematics, Gda\'nsk University, Ul.\ Wita Stwosza 57, 80-952 Gda\'nsk, Poland} 
\email{wrosicki@math.univ.gda.pl} 
\author{A.~Zastrow} 
\address{Institute of Mathematics, Gda\'nsk University, Ul.\ Wita Stwosza 57, 80-952 Gda\'nsk, Poland} 
\email{mataz@math.univ.gda.pl} 
\author{M.~\v Zeljko} 
\address{Institute of Mathematics, Physics and Mechanics, University of Ljubljana, P.O.Box 2964, Ljubljana 1001, Slovenia} 
\email{matjaz.zeljko@fmf.uni-lj.si} 

\date{March 29, 2008} 

\keywords{Embeddability, 
polyhedron, cone, 
suspension, 
Peano continuum, 
planar graph, 
local planarity} 
\subjclass[2000]{Primary: 57Q35; Secondary: 54C25, 55S15, 57N35} 
\begin{abstract} 
The main result of this paper is that if $X$ is a Peano continuum
such that its $n$-th cone $C^n(X)$ embeds into $\RR^{n+2}$ then
$X$ embeds into $S^2$. This solves a problem proposed by W. Rosicki.
\end{abstract} 
\maketitle

\section{Introduction} 
The classical Lefschetz-N\"obeling-Pontryagin Embedding Theorem \cite{Engelking1} asserts that
every compact metric space $X$ of dimension $n$ embedds into
$\RR^{2n+1}$. 
We are interested in the 
relationship between the
embeddability of $X$ and embeddability of its Cartesian product $X\times I^n$ with a cube $I^n$ 
(resp.\ its cone $C(X)$, iterated cone $C^n(X) = C(\ldots(C(X))\ldots)$, suspension $\Sigma(X)$).
Clearly, if $X$ embeds in $\RR^m$, then $X \times I^n$ and
$C^n(X)$ embed into $\RR^{n+m}$. However, sometimes they embed into lower-dimensional Euclidean space. Such is the case for the spheres
$S^n$, where $S^n$, $C(S^n) \cong B^{n+1}$ and $S^n \times I$ all embed into $\RR^{n+1}$.

Let $X$ be a Peano continuum. It was proved in \cite{Rosicki2} 
that if the cone $C(X)$ of $X$ embeds into $\RR^3$, then $X$ embeds into $S^2$. 
As a consequence, if the suspension $\Sigma(X)$ of $X$ embeds into $\RR^3$, then $X$ is planar. 
Note that for each $n\geq3$, there exists a Peano continuum $X_n$ such that $X_n$ is not embeddable in $S^n$,
whereas the cone $C(X_n)$ of $X_n$ is embeddable in $\RR^{n+1}$ (see \cite{Rosicki2}).

The main result of this paper is Theorem \ref{MainTheo} which solves a problem from \cite{Rosicki2}. Our proof is
based on the methods of \cite{Cauty1} and \cite{Rosicki2}.

\begin{figure} [htb] 
\begin{center} 
{\epsfxsize=46mm\epsfclipon\epsfbox{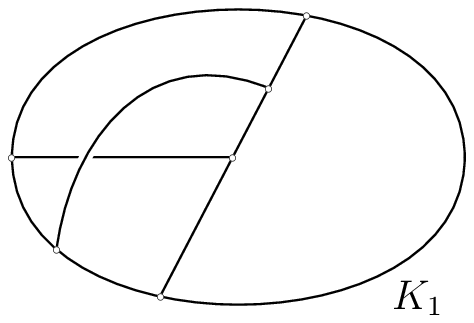}\epsfclipoff} \qquad
{\epsfxsize=46mm\epsfclipon\epsfbox{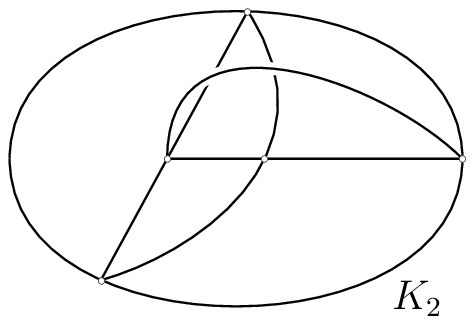}\epsfclipoff} \\
{\epsfxsize=46mm\epsfclipon\epsfbox{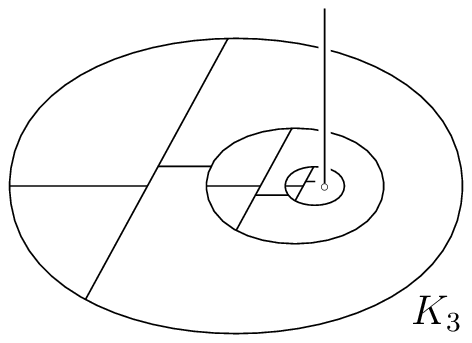}\epsfclipoff} \qquad
{\epsfxsize=46mm\epsfclipon\epsfbox{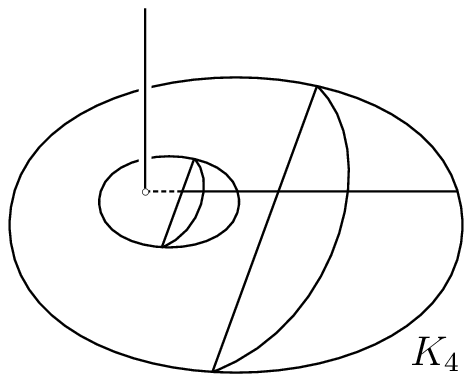}\epsfclipoff} 
\end{center} 
\caption{Kuratowski curves $K_1$, $K_2$, $K_3$, $K_4$} \label{Fig1}
\end{figure}

\begin{theo} \label{MainTheo}
Let $X$ be a Peano continuum. 
Suppose that for some $n\in\NN$, $C^n(X)$ is embeddable in $\RR^{n+2}$.
Then $X$ is embeddable in $S^2$.
\end{theo} 
Let $X$ be a Peano continuum. Claytor  \cite{Claytor2} proved
that $X$ is 
embeddable in $S^2$ if and only if $X$ does not contain any of the Kuratowski curves 
$K_1$, $K_2$, $K_3$, $K_4$ (see Figure \ref{Fig1}).

\section{Preliminaries}

A space $X$ is said to be \emph{planar} if $X$ is embeddable in $\RR^2$. We say that $X$ is \emph{locally planar} if for every point $x\in X$ there exists a neighbourhood $U_x$ of $x$ in $X$ such that $U_x$ is embeddable in $\RR^2$. 
Rosicki \cite[Theorem 1.1]{Rosicki1} proved  
that if a Peano continuum $X$ is embeddable in $\RR^3$ and $X$ is a nontrivial Cartesian product $X = Y\times Z$ then one of the factors is either an arc or a simple closed curve. 

Rosicki \cite{Rosicki1} also proved  that if a Peano continuum $X$ is embeddable in $\RR^3$ and is homeomorphic to the product $Y\times S^1$ then the factor $Y$ must be planar. Alternatively, if $X=Y\times[0,1]$ is embeddable in $\RR^3$ and $\check{H} ^1(X)=\check{H} ^2(X)=0$ then $Y$ must be planar.
Cauty \cite{Cauty1}, generalizing Rosicki \cite{Rosicki1},  
proved that for every $n>3$ and every Peano continuum $X$ such that $X\times I^{n-2}$ is embeddable into an
$n$-manifold, it follows that $X$ must be locally planar.
This theorem was stated earlier by Stubblefield \cite{Stubblefield1}. However, Burgess \cite{Burgess1} found a mistake in his proof.

Borsuk  \cite{Borsuk1} constructed an example of a locally connected, locally planar continuum $X$ which is not embeddable
into any surface. This continuum contains a sequence $(X_n)$ of subsets homeomorphic to Kuratowski curve $K_1$ which converge to an arc.
Cauty \cite{Cauty1} 
proved that $X\times I^{n-2}$ is not embeddable into any $n$-manifold so the converse to his
theorem does not hold.

\section{Local separation} 

We say that a subset $D\subset\RR^n$ \emph{locally separates}  $\RR^n$ at the point $x_0\in D$ into $k\in \NN$ components if there exists $\varepsilon>0$
such that for all $0<\delta<\varepsilon$, the set $B(x_0,\delta)\setminus D$ has exactly $k$ components $A_1$, \ldots, $A_k$ for which $x_0\in \overline{A_i}$, for
all $i\in\{1,\ldots, k\}$.

It is easy to prove the following lemma using similar methods as in the proof of Lemma 1 in \cite{Rosicki2}.

\begin{lemm} \label{L1} 
A homeomorphic image of any $n$-disk locally separates $\RR^{n+1}$ at any point of its interior into two components.
\end{lemm} 

Note that $C^n(X)=\sigma^{n-1} \ast X = \{xt+y(1-t);\ x\in \sigma^{n-1}, y\in X, t\in[0,1]\}$,
where $\sigma^{n-1}$ is an $(n-1)$-simplex.
Then   $\sigma^{n-1} \ast \{x\}$
is an $n$-ball
and
$\sigma^{n-1} \ast I$ is an $(n+1)$-ball. 
We consider $\sigma^{n-1}$ as a subset of $\sigma^{n-1}\ast X.$ 

\begin{lemm} \label{L2} 
Let $I_i$, $i\in\{1,\ldots, k\}$, $k>1$ be arcs with common endpoints and pairwise disjoint interiors and
$C_k=C^n(\bigcup_{i=1} ^k I_i)=\sigma^{n-1} \ast (\bigcup_{i=1} ^k I_i)$. Let 
$h\colon C_k\to\RR^{n+2}$ be
an embedding. Then $h(C_k)$ locally separates $\RR^{n+2}$ at any point $h(x_0)$, where $x_0$ is an interior
point of $\sigma^{n-1}$, into $k$ components
(where $\sigma^{n-1}$ is considered as a subset of $C_k$).
\end{lemm} 

{\sl Proof.} 
The proof is by induction on $k$. 
 If $k=2$, then 
 $C_2=\sigma^{n-1}\ast S^0 \ast S^0$ hence
 $h(C_2)$ locally separates $\RR^{n+2}$ at $h(x_0)$ into
two components, by Lemma \ref{L1}.

Assume that Lemma \ref{L2} holds for $k-1$. Choose $\varepsilon>0$ smaller than the distance between $h(x_0)$ and the 
image of $\partial\sigma^{n-1}\ast (\bigcup_{i=1} ^k I_i)$.
Let $\delta>0$ be so small that 
$$D_k=h(C_{k}\cap B(x_0,\delta))\subset B(h(x_0),\varepsilon).$$

There exists an open connected set $U_k\subset\RR^{n+2}$ 
such that 
$D_k=U_k\cap h(C_k)$.
Consider the exact sequence of the pair $(U_k,U_k\setminus D_k)$:
$$\to H_1(U_k)\to H_1(U_k,U_k\setminus D_k)\to H_0(U_k\setminus D_k)\to H_0(U_k) \to H_0(U_k,U_k\setminus D_k)\to 0.$$
Since $U_k$ is an open $(n+2)$-manifold, $H_1(U_k)\cong \check H_c^{n+1}(U_k)$ by the
Poincar\'e duality,
where $\check H_c$ denotes the \v Cech cohomology with compact supports. 
Also $H_1(U_k,U_k\setminus D_k)\cong \check H_c^{n+1}(D_k)$ (see \cite[VIII, 7.14]{Dold1}, where $L=\emptyset$, $K=D_k$ and $X=U_k$).

We know that $H_0(U_k,U_k\setminus D_k)=0$ because $U_k$ is arc-connected and $U_k\setminus D_k\ne\emptyset$.
Therefore we can consider the exact sequence
$$\to \check H_c^{n+1}(U_k)\to \check H_c^{n+1}(D_k)\to H_0(U_k\setminus D_k)\to H_0(U_k) \to 0.$$
Next we show by induction that the map 
$\check H_c^{n+1}(U_k)\to \check H_c^{n+1}(D_k)$ is trivial. If $k=2$ then $D_k$ is an open $(n+1)$-ball. Then 
$H_0(U_k\setminus D_k)\cong \ZZ^2$, by Lemma \ref{L1}. Since $\check H_c^{n+1}(D_k)\cong\ZZ$ and $H_0(U_k)\cong\ZZ$,
we 
obtain the exact sequence 
$$\check H_c^{n+1}(U_k)\to \ZZ\to \ZZ^2\to \ZZ \to 0.$$
Hence the map $\check H_c^{n+1}(U_k)\to \check H_c^{n+1}(D_k)$
is indeed trivial, as asserted.

Since $\check H_c^{n+1}(D_2)\cong\ZZ$, we obtain by induction that 
$\check H_c^{n+1}(D_k)\cong \check H_c^{n+1}(D_{k-1})\oplus \check H_c^{n+1}(D_2')\cong \ZZ^{k-2}\oplus\ZZ$, 
where $D_2'=h(C^n(I_1\cup I_k) \cap B(x_0,\delta))$. 

The map
$\check H_c^{n+1}(U_k)\to  \check H_c^{n+1}(D_k)\cong  \check H_c^{n+1}(h(D_{k-1}))\oplus \check H_c^{n+1}(D_2')$ 
is trivial because both of its coordinates are trivial, by inductive hypothesis.

Therefore the sequence 
$$0\to \check H_c^{n+1}(D_k)\to H_0(U_k\setminus D_k)\to H_0(U_k)\to0$$
is exact. So the sequence 
$$0\to \ZZ^{k-1}\to H_0(U_k\setminus D_k)\to \ZZ\to0$$
is also exact. Hence
$H_0(U_k\setminus D_k)\cong\ZZ^k$  and $U_k\setminus D_k$ has $k$ components.

The point $h(x_0)$ belongs to the closure of each of them. Indeed, if $X_k$ is $D_k$ with a small open neighbourhood of $h(x_0)$ removed
then $\check H_c^{n+1}(X_k)\cong0$ and the sequence
$$0\to H_0(U_k\setminus X_k)\to H_0(U_k)\to0$$
is exact, therefore $H_0(U_k\setminus X_k)\cong\ZZ$.\qed

\section{Proof of Theorem 1.1}

We shall need two more lemmata:

\begin{lemm} \label{L3} 
Consider the Kuratowski curve $K_1$ and let $n\in\NN$. Then $C^n(K_1)$ is not embeddable in $\RR^{n+2}$.
\end{lemm} 

\begin{figure} [htb] 
\begin{center} 
{\epsfxsize=49mm\epsfclipon\epsfbox{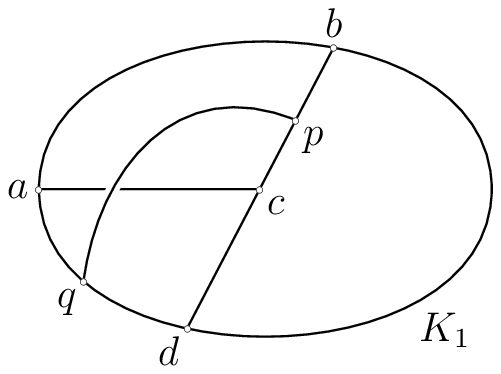}\epsfclipoff} \qquad\qquad
{\epsfxsize=50mm\epsfclipon\epsfbox{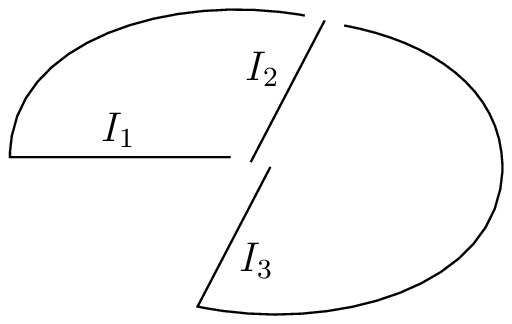}\epsfclipoff} 
\end{center} 
\caption{Kuratowski curve $K_1$} \label{Fig2}
\end{figure}

{\sl Proof.} 
Suppose to the contrary, that there exists an embedding $h\colon C^n(K_1)\to\RR^{n+2}$.
Consider $K_1\subset \RR^3$ and denote (see Figure \ref{Fig2})
$$ I_1 = [c,a] \cup [a,b], \  I_2 = [c,p] \cup [p,b], \mbox{ and } I_3 = [c,d] \cup [d,b].$$


\begin{figure} [htb] 
\begin{center} 
{\epsfxsize=46mm\epsfclipon\epsfbox{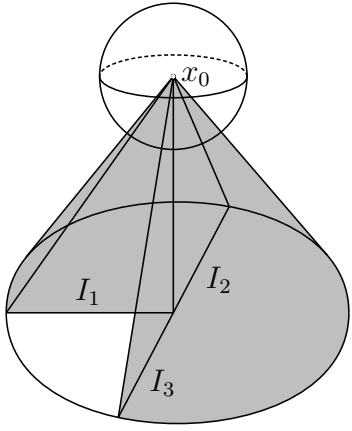}\epsfclipoff} 
\end{center} 
\caption{Local separation at $h(x_0)$} \label{Fig3}
\end{figure}

If $X=\bigcup_i I_i$, then $\sigma^{n-1} \ast X=\bigcup_i (\sigma^{n-1}\ast I_i)$ is a union of $(n+1)$-disks.
Let $x_0\in\Int \sigma^{n-1}$ and choose $\varepsilon>0$ so that (see Figure \ref{Fig3})
\begin{eqnarray*} 
C_1 & = & h(\sigma^{n-1} \ast ( I_1 \cup I_3 )) \mbox{ locally separates $B(h(x_0),\varepsilon)$ into } B_1, A_1 \mbox{ at }h(x_0), \\
C_2 & = & h(\sigma^{n-1} \ast ( I_1 \cup I_2 )) \mbox{ locally separates $B(h(x_0),\varepsilon)$ into } B_2, A_2 \mbox{ at }h(x_0), \\
C_3 & = & h(\sigma^{n-1} \ast ( I_2 \cup I_3 )) \mbox{ locally separates $B(h(x_0),\varepsilon)$ into } B_3, A_3 \mbox{ at }h(x_0). 
\end{eqnarray*}


By Lemma \ref{L2} we have that $C=h(C^n(I_1\cup I_2\cup I_3))=h(\sigma^{n-1} \ast \bigcup_{i=1} ^3 I_i)=h(\bigcup_{i=1} ^3 \sigma^{n-1} \ast I_i)$
locally separates $B(h(x_0),\varepsilon)$ into three components. We will show that we can adopt the notation for these three 
components to be $B_1$, $A_2$ and $A_3$.

We use abstract linear combinations for describing our joins, e.g.
$$\sigma^{n-1} \ast K_1 = \{xt+y(1-t);\ x\in \sigma^{n-1}, y\in K_1, t\in[0,1]\}.$$

For $\sigma^{n-1} \subset \sigma^{n-1} \ast K_1$, we have that $h(\sigma^{n-1})$ is a
subset of $C_1$, but that $h|_{\sigma^{n-1} \ast I_2}$ maps all linear combinations with $t \neq 1$, but sufficiently close
to $1$,  to a subset that
is 
connected but disjoint from $C_1$. Hence this subset can only
be contained either in
$A_1$ or in $B_1$. We may assume that it is in $A_1$.
Since the entire
neighbourhood of $\sigma^{n-1}$ in $\sigma^{n-1} \ast I_2$ is mapped by $h$ into $A_1$, we have
$h(\sigma^{n-1} \ast I_2 )\cap B_1=\emptyset$, provided $\varepsilon>0$ is small enough.
Then $B_1$ is not divided by $C$, so it is one of the three components.

Analogously, by considering $C_2$ (resp.~$C_3$) we can make sure that $A_2$ and $A_3$ are the other two components and that
$h(\sigma^{n-1} \ast I_3 )\cap A_2=\emptyset$
and $h(\sigma^{n-1} \ast I_1 )\cap A_3=\emptyset$.
Since $C\cup B_1\cup A_2\cup A_3$ and  $C\cup A_1\cup B_1$ are both disjoint decompositions of a neighbourhood of $h(x_0)$, 
the set $h(\sigma^{n-1} \ast I_2 )\cup C_1$ separates the component $A_1$ into components $A_2$ and $A_3$.

Note that
$$x_0\ast K_1 = \{x_0t+x(1-t);\ x \in K_1, t\in[0,1]\} \subset C^n(K_1).$$
Choose $t_0$ near 1 so that
$$h(\{x_0t+x(1-t);\ x \in K_1, t\geq t_0\} )\subset B(h(x_0),\varepsilon).$$
Let $p'=h(x_0t_0+p(1-t_0))\in A_1$. 
The arc 
$H=h(\{x_0t_0+x(1-t_0);\ x \in (p,q)\} )$ is contained in $B(h(x_0),\varepsilon)\setminus h(C)$.
Therefore
points $p'$ and $q'=h(x_0t_0+q(1-t_0))$ are in the same component. 
Hence $q'\in A_2$ or $q'\in A_3$.
So the arc $I=h(\{x_0t_0+x(1-t_0);\ x \in (a,q]\cup[q,d)\} )$ is contained either in $A_2$ or in $A_3$.
But this yields a contradiction since $a'=h(x_0t_0+a(1-t_0))\notin \overline{A_3}$ (so $I\not \subset A_3$)
and $d'=h(x_0t_0+a(1-t_0))\notin \overline{A_2}$ (so $I\not \subset A_2$). \qed

The proof of the next lemma can be obtained by changing the proof of \cite[Lemma 4]{Rosicki2} 
in the same way as we did it for the proof of Lemma 2.3 using the proof of \cite[Lemma 3]{Rosicki2}.

\begin{lemm} \label{L4} 
Consider the Kuratowski curve $K_2$ and let $n\in\NN$. Then 
$C^n(K_2)$ is not embeddable in $\RR^{n+2}$.
\end{lemm}

\bigskip

{\sl Proof of Theorem \ref{MainTheo}.} 
By Claytor's theorem (see \cite{Claytor1}, \cite{Claytor2}),
it suffices to show that $C^n(K_i)$ 
is not embeddable into $\RR^{n+2}$ for any 
$i\in\{1,2,3,4\}$. 
Now, Cauty \cite{Cauty1} proved that $K_i\times I^n$
is not embeddable into $\RR^{n+2}$ for any
$i\in\{3,4\}$. 
Therefore also $C^n(K_i)$ is not embeddable into $\RR^{n+2}$ for any
$i\in\{3,4\}$.
Hence we only have to consider the cases $i=1$ and $i=2$.
The proof is now completed by application of Lemmata 4.1 and 4.2.\qed 

\section{Epilogue} 

Repov\v s, Skopenkov and \v S\v cepin  \cite{ReSkSh1} proved that if
$X\times I$ PL embeds into $\RR^{n+1}$, where $X$ is either
an acyclic polyhedron and $\dim X\leq \frac{2n} {3} -1$ or
a homologically $(2\dim X-n-1)$-connected manifold and $\dim X\leq \frac{2n} {3} -1$ or
a collapsible polyhedron, then $X$ PL embeds into $\RR^n$.

\begin{ques}
What can one say about embeddability of $X$ into Euclidean spaces
if one considers $C(X)$ or $C^n(X)$ or $\Sigma(X)$ or $\Sigma^n(X)$ instead of $X\times I$ for $X$ in \cite{ReSkSh1}?
 \end{ques}

It follows by \cite{ReSkSh1} that if $X$ is a contractible polyhedron such that $X\times I$ embeds into $\RR^{n+1}$ then $X$ embeds into
$\RR^n$. So if $X$ is contractible and $C(X)\subset \RR^{n+1}$ then $X$ embeds into $\RR^n$. 

Note that there exists a
polyhedron $P_n$ such that 
$P_n$ is not embeddable into $\RR^n$ but $C^2(P_n)$ is embeddable in $\RR^{n+2}$.
Namely, Cannon \cite{Cannon1} proved that if $H^n$
is a homology $n$-sphere then its double suspension $\Sigma^{2}(H^n)$is the $(n+2)$-sphere (see \cite{Daverman1} \cite{Edwards1} for a far reaching generalization of this result). So
if $P_n=H^n\setminus B^n$ where $B^n$ is an $n$-ball then the double cone
$C^2(P_n)$ embeds in $\RR^{n+2}$. 
The polyhedron $P_n$ is acyclic but not contractible.

\begin{ques}
Does there exist a contractible $n$-dimensional polyhedron $X^n$ such that $C^k(X^n)$ embeds into
$\RR^{n+k}$, but $X^n$ does not embed into $\RR^n$?\end{ques}

In \cite[Theorem 2]{Rosicki2} contractible continua $X_n$ were constructed,  such that $X_n$ is not embeddable in 
$\RR^n$, $C(X_n)$ is embeddable in $\RR^{n+1}$, and $X_n$ is  not a polyhedron.
By \cite{ReSkSh1}, if $X$ is an $n$-polyhedron then $X\times I$ embeds into $\RR^{2n+1}$.
If $X$ is an
$n$-polyhedron then $C(X)$ need not
embed into $\RR^{2n+1}$.
For example, the Kuratowski curves
$K_1$ and $K_2$ are $1$-polyhedra but the cones
$C(K_1)$ and $C(K_2)$ do not embed into $\RR^3$.

\begin{ques}
Suppose that $X$ is a compact contractible $n$-dimensional polyhedron. 
Does the cone $C(X)$ embed into $\RR^{2n+1}$? Does the same hold if $X$ is only acyclic?
\end{ques}

\section{Acknowledgements}
This research was supported by 
Polish--Slovenian  grant BI-PL 12/2004-2005, 
Polish grant N20100831/0524, 
SRA program P1-0292-0101-04,
and SRA project J1-9643-0101.
The authors acknowledge R.~Cauty \cite{Cauty2}
for the hint communicated to the second author, and the referee for comments and suggestions. 


\end{document}